\documentclass[12pt,a4paper]{article}
\usepackage[T1]{fontenc}
\usepackage{mathptmx}
\usepackage{courier}
\usepackage[dvips]{graphicx}
\usepackage{rotating}
\usepackage{psfrag}
\usepackage{amsmath,amssymb}        
\usepackage{latexsym}
\usepackage{url}

\title{A Non-Iterative Transformation Method\\ Applied to Boundary-Layer Flows\\ of Non-Newtonian Fluids Past a Flat Plate}
\author{Riccardo Fazio\footnote{e-mail: rfazio@unime.it \ \ \ home-page: http://mat521.unime.it/fazio} \\
ORCID: 0000-0003-0825-0162 \\
Department of Mathematics, Computer Science\\ Physical Sciences and Earth Sciences (MIFT),\\
University of Messina \\
Viale F. Stagno D'Alcontres, 31 \\
98166 Messina, Italy}
\date{August 21, 2020}
\begin{document}
\maketitle
\begin{abstract}
In this paper, we define a non-iterative transformation method for boundary-layer flows of non-Newtonian fluids past a flat plate.
The problem to be solved is an extended Blasius problem depending on a parameter.
This method allows us to solve numerically the extended Blasius problem by solving a related initial value problem and then rescaling the obtained numerical solution. 
Therefore, it is a non-iterative initial value method.
We find that our computed numerical results, for a wide range of the parameter involved, are in very good agreement with the data reported in the literature.
\end{abstract}
\bigskip

\bigskip


\section{Introduction.}
The most valuable development in fluid mechanics within the 20th century was, certainly, the concept of boundary layer flows introduced by Prandtl \cite{Prandtl:1904:UFK}.
A boundary layer is that layer of fluid that forms near a surface that is subject to the fluid flow.
Therefore, boundary layers occur in several practical problems like flows: inside the blood vessels, attached to aeroplane wings, in irrigation channels, near earth's surface and around buildings due to winds, in sewer pipes, around a moving car and many others.  
There are two simplest problems in this contest in the literature. 
The first one describes the flow along with a horizontal flat motionless plate due to a constant free stream, see Blasius \cite{Blasius:1908:GFK}. 
The second is flow induced by a horizontal flat plate moving with constant velocity inside a quiet fluid, see Sakiadis \cite{Sakiadis:1961:BLBa,Sakiadis:1961:BLBb}.
In the first problem the fluid velocity increases from zero at the plate, no-slip boundary condition, to the mainstream velocity far away from the plate.
In the second problem, the fluid velocity is equal to the plate velocity at the plate, no-slip condition, and decreases to zero far away from the plate.  
In both cases, the practical interest is to compute the shear at the plate (skin friction), this datum is essential to the determination of the viscous drag on the plate, see Schlichting \cite{Schlichting:2000:BLT}.
It turns out that the increase of the wall shear in Sakiadis solution with respect to Blasius solution is about 33.64\%, see for instance \cite{Fazio:2015:ITM}.
In the above contest the Blasius boundary value problem (BVP) is given by
\begin{align}\label{eq:Blasius} 
& {\displaystyle \frac{d^3 f}{d \eta^3}} + \frac{1}{2} \; f
{\displaystyle \frac{d^{2}f}{d\eta^2}} = 0 \nonumber \\[-1ex]
&\\[-1ex]
& f(0) = {\displaystyle \frac{df}{d\eta}}(0) = 0, \qquad
{\displaystyle \frac{df}{d\eta}}(\eta) \rightarrow 1 \quad \mbox{as}
\quad \eta \rightarrow \infty \ . \nonumber 
\end{align}
This is a BVP defined on the semi-infinite real axis $[0, \infty)$.
Weyl \cite{Weyl:1942:DES} proved that the unique solution of (\ref{eq:Blasius}) has a monotone decreasing positive second order derivative on all $ [0, \infty) $ and approaches to zero as $ \eta $ goes to infinity.
In (\ref{eq:Blasius}) the governing differential equation and the two boundary conditions
at the origin are invariant with respect to the scaling group of point transformations
\begin{eqnarray}
& \eta^* = \lambda^{-\alpha}, \qquad f^* = \lambda^{\alpha} f
\label{eq:scaling:Blasius}
\end{eqnarray}
where $ \alpha $ is a nonzero constant, while the asymptotic boundary condition is not invariant. 
The interest related to this invariance properties is both analytical and numerical. 
As a consequence to that transformation, a simple existence and uniqueness Theorem
was given by J. Serrin, this is cited, for instance, by Meyer \cite[pp. 104-105]{Meyer:1971:IMF}. 
AS far as a numerical viewpoint is concerned, a non-iterative transformation method (ITM)
reducing the solution of (\ref{eq:Blasius}) to the solution of a related initial value problem (IVP) was deduced by T\"opfer \cite{Topfer:1912:BAB}. 
Moreover, the mentioned invariance properties are basic to the error analysis of the truncated boundary solution proposed by Rubel \cite{Rubel:1955:EET}, see Fazio \cite{Fazio:2002:SFB} for the full details.
Blasius problem was used, recently, by Boyd \cite{Boyd:2008:BFC} as an example were some preliminary analysis allowed researchers of the past to solve problems, governed by partial differential equations, that, before the computer invention, might be otherwise impossible to face.

We have applied a non-ITM to several problems of practical interest. 
For istance, to the Blasius equation with slip boundary condition, arising within the study of gas and liquid flows at the micro-scale regime \cite{Gad-el-Hak:1999:FMM,Martin:2001:BBL}, as reported in \cite{Fazio:2009:NTM}.
To the Blasius equation with moving wall considered by Ishak et al. \cite{Ishak:2007:BLM} or surface gasification studied by Emmons \cite{Emmons:1956:FCL} and recently by Lu and Law \cite{Lu:2014:ISB} or slip boundary conditions investigated by Gad-el-Hak \cite{Gad-el-Hak:1999:FMM} or Martin and Boyd \cite{Martin:2001:BBL}, see Fazio \cite{Fazio:2016:NIT} for details.
In this contest, we found a way to solve non-iteratively the Sakiadis problem \cite{Sakiadis:1961:BLBa,Sakiadis:1961:BLBb}.
A recent review dealing with the derivation and application of non-ITM can be be found, by the interested reader, in \cite{Fazio:2019:NIT}.

Moreover, T\"opfer's method has been extended to classes of problems in boundary layer theory involving one, or more than one, physical parameter.
Na \cite{Na:1970:IVM} was the first to study such an extension, see also NA \cite[Chapters 8-9]{Na:1979:CME} for an extensive survey of this subject.

Finally, an iterative extension of the transformation method has been defined, for the numerical solution of free BVPs, by Fazio and Evans \cite{Fazio:1990:SNA} see also the applications studied by Fazio \cite{Fazio:1990:NVT,Fazio:1991:ITM,Fazio:1992:ITM,Fazio:1997:NTE,Fazio:2001:ITM}. 
This iterative extension was used to solve numerically several problems of practical interest: 
free boundary problems \cite{Fazio:1990:SNA,Fazio:1997:NTE,Fazio:1998:SAN},
a moving boundary hyperbolic problem \cite{Fazio:1992:MBH},
Homann and Hiemenz's problems governed by the Falkner-Skan equation in \cite{Fazio:1994:FSE},
one-dimensional parabolic moving boundary problems \cite{Fazio:2001:ITM}, two variants of the Blasius problem \cite{Fazio:2009:NTM}, namely: a boundary layer problem over moving surfaces, studied first by Klemp and Acrivos \cite{Klemp:1972:MBL}, and a boundary layer problem with slip boundary condition, that can be found in the study of gas and liquid flows at the micro-scale regime \cite{Gad-el-Hak:1999:FMM,Martin:2001:BBL}, parabolic problems defined on unbounded domains \cite{Fazio:2010:MBF} and, recently, see \cite{Fazio:2015:ITM}, an additional variant of the Blasius problem in boundary layer theory: the so-called Sakiadis problem \cite{Sakiadis:1961:BLBa,Sakiadis:1961:BLBb}.
As far as the ITM is concerned, a recent review dealing with all the cited problems can be be found in \cite{Fazio:2019:ITM}.

\section{Boundary layer of non-Newtonian fluids past a flat plate}
The boundary layer of non-Newtonian fluids flowing past a flat plate is described, using suitable similarity variables, by an extended Blasius problem 
\begin{align}\label{eq:ExBlasius} 
& P (P+1) \; {\displaystyle \frac{d^3 f}{d \eta^3}} + f
{\displaystyle \frac{d^{2}f}{d\eta^2}}^{(2-P)} = 0 \nonumber \\[-1ex]
&\\[-1ex]
& f(0) = {\displaystyle \frac{df}{d\eta}}(0) = 0, \qquad
{\displaystyle \frac{df}{d\eta}}(\eta) \rightarrow 1 \quad \mbox{as}
\quad \eta \rightarrow \infty \ , \nonumber 
\end{align}
where $P$ verifies the one-side conditions $0 \le P $, see Acrivos et al. \cite{Acrivos:1960:MHT}.
We remark that our problem (\ref{eq:ExBlasius}) reduces, for $P = 1$, to the celebrated Blasius problem (\ref{eq:Blasius}).

\subsection{The non-ITM} 
Let us remark here that, a non-ITM can be applied to the Blasius problem (\ref{eq:Blasius}) because the governing differential equation and the two boundary conditions at $\eta = 0$ are invariant with respect to and the asymptotic boundary condition is not invariant under the scaling transformation (\ref{eq:scaling:Blasius}).
In order to apply a non-ITM to the BVP (\ref{eq:ExBlasius}) we investigate its invariance under the scaling group of point transformation
\begin{equation}\label{eq:scaling}
f^* = \lambda f \ , \qquad \eta^* = \lambda^{\delta} \eta \ .   
\end{equation}
We find that the extended Blasius problem (\ref{eq:ExBlasius}) is invariant with respect to (\ref{eq:scaling}) iff
\begin{equation}\label{eq:scaling:condition}
\delta = \frac{P-2}{2 P-1} \ .
\end{equation}
Now, we can integrate the extended Blasius equation in (\ref{eq:ExBlasius}) written in the star variables on $[0, \eta^*_\infty]$, where $\eta^*_\infty$ is a truncated boundary that can be found by trial, with initial conditions
\begin{equation}\label{eq:ICs2}
f^*(0) = \frac{df^*}{d\eta^*}(0) = 0 \ , \quad \frac{d^2f^*}{d\eta^{*2}}(0) = 1 \ ,
\end{equation}
in order to compute an approximation $\frac{df^*}{d\eta^*}(\eta^*_\infty)$ for $\frac{df^*}{d\eta^*}(\infty)$ and the corresponding value of $\lambda$ according to the equation
\begin{equation}\label{eq:lambda}
\lambda = \left[ \frac{d f^*}{d \eta^{*}}(\eta^*_\infty) \right]^{1/(\delta-1} \ .   
\end{equation} 
Computed the value of $\lambda$ by equation (\ref{eq:lambda}), we can find the missed initial condition by the equation
\begin{equation}\label{eq:MIC}
\frac{d^2f}{d\eta^{2}}(0) =  \lambda^{2\delta-1}\frac{d^2f^*}{d\eta^{*2}}(0) \ .
\end{equation}
Finally, we can compute the numerical solution of the original BVP (\ref{eq:ExBlasius}) by rescaling the numerical solution of the IVP.
Therefore, so doing we find the solution of a given BVP by solving a related IVP.

\section{Numerical results}
To compute the numerical solution, we used an explicit eight order Runge-Kutta method \cite[p. 180]{Butcher:NMO:2003} with constant step size.
In table \ref{tab:missingIC} we report the missing initial conditions $\frac{d^2f}{d\eta^2}(0)$ for several values of the parameter $P$.
\begin{table}[!hbt]
\caption{Missing initial conditions for different values of $P$.}
\vspace{.5cm}
\renewcommand\arraystretch{1.3}
	\centering
		\begin{tabular}{r@{.}lr@{.}lr@{.}lr@{.}l}
\hline 
\multicolumn{2}{c}%
{$P$} &  
\multicolumn{2}{c}%
{$\frac{d^2f}{d\eta^2}(0)$} by \cite{Acrivos:1960:MHT} &  
\multicolumn{2}{c}%
{$\frac{d^2f}{d\eta^2}(0)$} from (\ref{eq:Pohlhausen}) &  
\multicolumn{2}{c}%
{$\frac{d^2f}{d\eta^2}(0)$} \\[1.2ex]
\hline
0 & 05 & 1 & 400938 & 0 & 214892 & 1 & 540752 \\
0 & 1  & 0 & 729857 & 0 & 221302 & 0 & 826478 \\  
0 & 2  & 0 & 505623 & 0 & 237305 & 0 & 490342 \\
0 & 3  & 0 & 354290 & 0 & 244046 & 0 & 391515 \\
0 & 4  &   &        &   &        & 0 & 350396 \\
0 & 5  & 0 & 331200 & 0 & 268324 &   &   \\
0 & 6  &   &        &   &        & 0 & 3239457 \\
0 & 7  &   &        &   &        & 0 & 3220337 \\
0 & 8  &   &        &   &        & 0 & 323544 \\
0 & 9  &   &        &   &        & 0 & 327139 \\
1 &    & 0 & 33206  & 0 & 323    & 0 & 332057  \\
1 & 5  & 0 & 363215 & 0 & 384047 & 0 & 398432 \\
\hline			
		\end{tabular}
	\label{tab:missingIC}
\end{table}
In table \ref{tab:missingIC} we compare our results with those reported in the paper by Acrivos et al. \cite{Acrivos:1960:MHT} and with those computed by the Pohlhausen method using the formula
\begin{equation}\label{eq:Pohlhausen}
{\displaystyle \frac{d^2f}{d\eta^2}(0)} = \left[\frac{39}{280} \frac{1.5}{P+1}\right]^{\frac{P^2}{P+1}} \ .
\end{equation}
From these results we can realize that the Pohlhausen method is inaccurate for small values of $P$ and that our numerical results are in very good agreement with those found by Acrivos et al.  \cite{Acrivos:1960:MHT}.
As remarked by Acrivos et al. \cite{Acrivos:1960:MHT}, for the range $P > 2$ the boundary layer flow is not an asymptotic state of laminar motion which is approached as the mainstream velocity is made sufficiently large. 
Therefore, we have neglected to investigate the range of $P \ge 2$.

\begin{figure}[!hbt]
	\centering
\psfragscanon 
\psfrag{e}[1][]{$\eta, \eta^*$}  
\psfrag{c}[1][]{$\frac{df^*}{d\eta^*}$,$\frac{d^2f^*}{d{\eta^*}^2}$}  
\psfrag{u}[1][]{$\frac{df}{d\eta}$,$\frac{d^2f}{d\eta^2}$}  
\includegraphics[width=14cm,height=8.5cm]{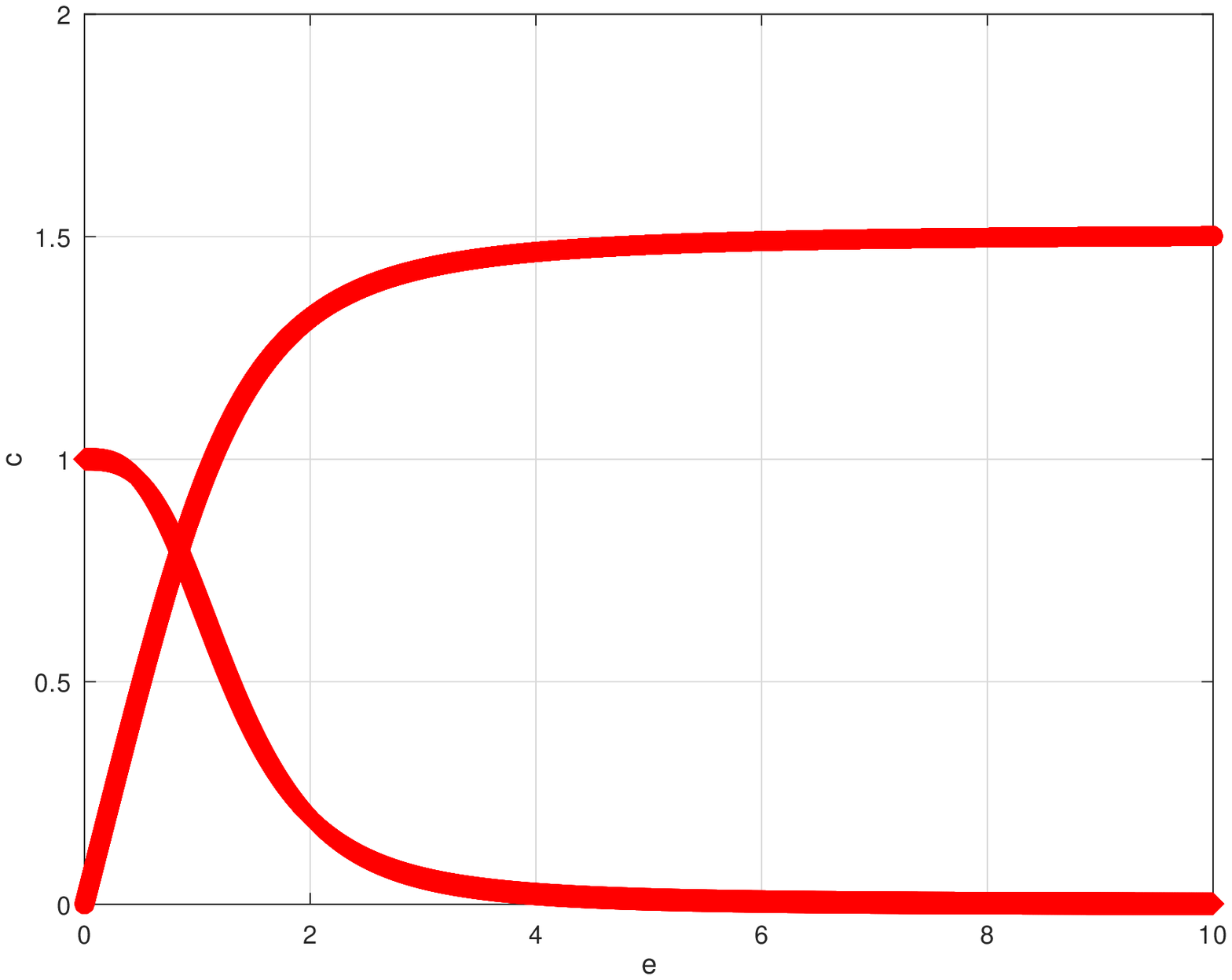} \\
\includegraphics[width=14cm,height=8.5cm]{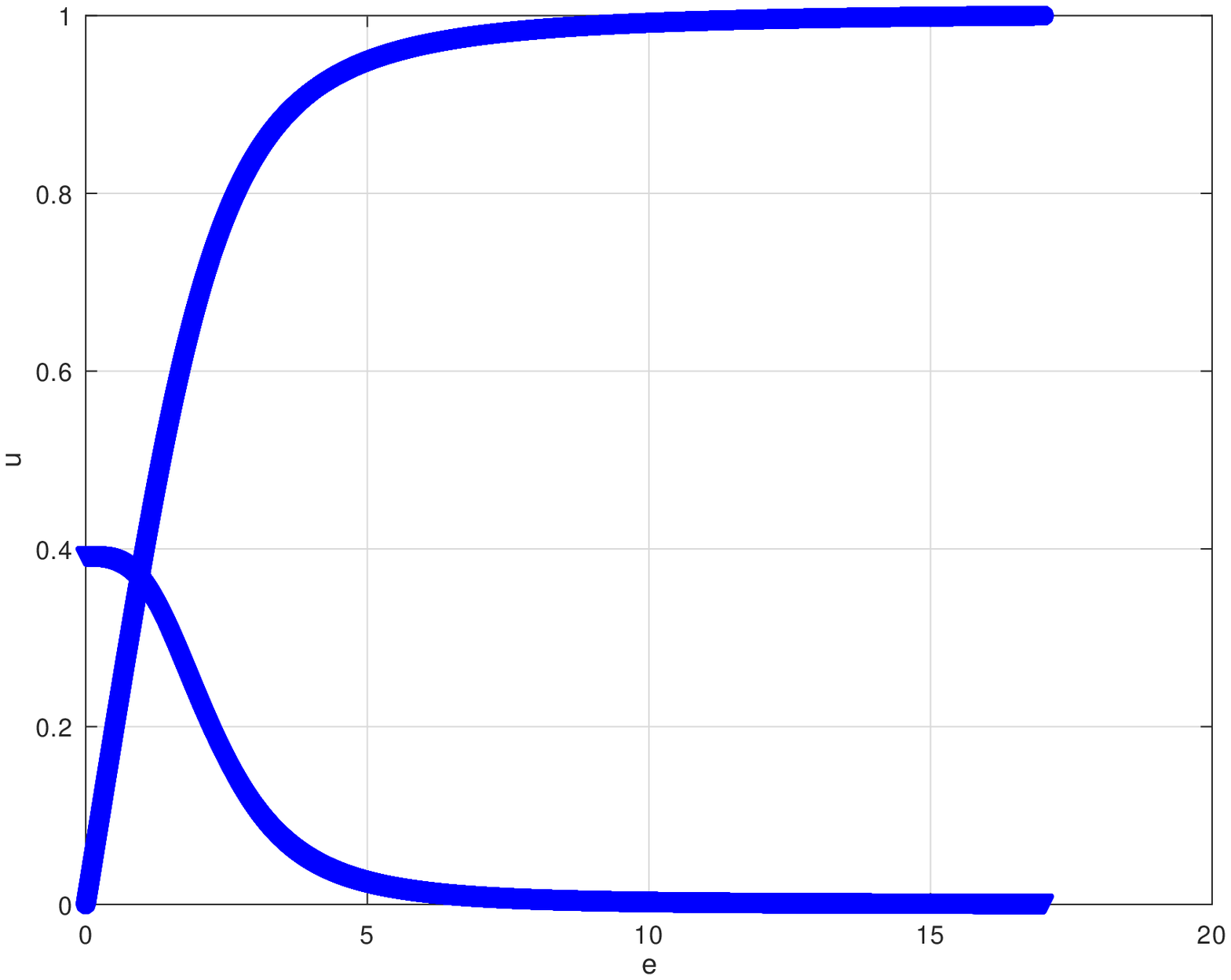}
\caption{Numerical results computed with the non-ITM for (\ref{eq:ExBlasius}) with $P =0.3$. The starred variables problem, top frame, and, bottom frame, the original problem solution components found after rescaling.}
	\label{fig:nBlasius}
\end{figure}

Figure \ref{fig:nBlasius} shows the solution of the extended Blasius problem, describing the behaviour of a boundary layer flow due to a moving flat surface immersed in an otherwise quiescent fluid, corresponding to $P=0.3$.
For the results shown in this figure we used an eight order Runge-Kutta method \cite[p. 180]{Butcher:NMO:2003} with constan step size $\Delta \eta = 0.001$ and a truncated boundary $\eta^*_{\infty} = 10$.
From figure\ref{fig:nBlasius} it is evident that $\eta^*_{\infty} < \eta_{\infty}$ and this is, of course, convenient from a numerical viewpoint because it means that we make let computational effort to computed the numerical solution of the original problem. 

The case $P=1$, as mentioned before, is the Blasius problem (\ref{eq:Blasius}).
In this case, our non-ITM reduces to the original method developed by T\"opfer \cite{Topfer:1912:BAB}, and the computed skin friction coefficient value, namely $0.332057336215$, obtained with $\Delta \eta = 0.001$ and $\eta^*_{\infty} = 10$, is in very good agreement with the values available in the literature, see for instance the value $0.332057336215$ computed by Fazio \cite{Fazio:1992:BPF} by a free boundary formulation of the Blasius problem or the value $0.33205733621519630$ computed by Boyd \cite{Boyd:1999:BFC} who believes that all the decimal digits to be correct.

\section{Final remarks and conclusions.}
The main result of this paper is the extension of the non-ITM, proposed by T\"opfer \cite{Topfer:1912:BAB} and defined for the numerical solution of the well-known Blasius problem \cite{Blasius:1908:GFK}, to an extended Blasius problem.
This method allows us to solve numerically the extended Blasius problem by solving a related initial value problem and then rescaling the obtained numerical solution. 
The obtained numerical results were computed by an eight order Runge-Kutta method \cite[p. 180]{Butcher:NMO:2003} to avoid the necessity, for accuracy reasons, to be forced to use very small step sizes.

\vspace{1.5cm}

\noindent {\bf Acknowledgement.} {This research was 
partially supported by the University of Messina and by the GNCS of INDAM.}


\end{document}